\documentclass{amsart}
\usepackage{latexsym,amssymb,amsthm,amsmath,verbatim}
\usepackage{lscape,pdflscape}
\usepackage{enumerate}
\usepackage{tikz}
\usetikzlibrary{arrows,positioning,chains,matrix,scopes}

\newtheorem{thm}{Theorem}[section]
\newtheorem{prop}[thm]{Proposition}
\newtheorem{cor}[thm]{Corollary}
\newtheorem{lem}[thm]{Lemma}

\theoremstyle{definition}

\newtheorem{example}[thm]{Example}

\begin{document}

%\cyh
%\dedicatory{Dedicated to Donald the Great}
\title[Relations between a topological game and the $G_{\delta}$-diagonal property]
{Relations between a topological game and the $G_{\delta}$-diagonal property}
\author[L. F. Aurichi]
{Leandro F. Aurichi$^1$}
\thanks{$^1$ Supported by FAPESP (2015/25725)}
\address []{Instituto de Ci\^encias Matem\'aticas e de Computa\c c\~ao,
Universidade de S\~ao Paulo, Caixa Postal 668,
S\~ao Carlos, SP, 13560-970, Brazil} \email{aurichi@icmc.usp.br}

\author[D. A. Lara]{Dione A. Lara$^2$}
\thanks{$^2$ Supported by CAPES (DS- 6111655/D)}
\address{Instituto de Ci\^encias Matem\'aticas e de Computa\c c\~ao,
Universidade de S\~ao Paulo, Caixa Postal 668,
S\~ao Carlos, SP, 13560-970, Brazil}
\email{dione@icmc.usp.br}

%\keywords{topological games, selection principles, productively countably tightness, Alster spaces, $G_\delta$-topology, bornology }

%\subjclass[2010]{Primary 54D20; Secondary 54G99, 54A10}

%\hjmtitle

\begin{abstract}
%For a cardinal $\kappa$ we present a selection principle that is equivalent to the space $X$ having the $G_{\kappa}$-diagonal property and we discuss some games related with that property. 

We present a selection principle $S_1(\mathcal{O},\mathcal{H})$ that characterizes the $G_{\delta}$-diagonal property. We also present a topological game induced by this selection principle and we study the relations between this game and the $G_{\delta}$-property. Finally, we give some applications and examples.

\end{abstract}
 
\maketitle

%%%Our Macros
\newcommand{\mP}{\mathcal{P}}
\newcommand{\mA}{\mathcal{A}}
\newcommand{\mC}{\mathcal{C}}
\newcommand{\mF}{\mathcal{F}}
\newcommand{\mB}{\mathcal{B}}
\newcommand{\mH}{\mathcal{H}}
\newcommand{\mU}{\mathcal{U}}
\newcommand{\mO}{\mathcal{O}}
\newcommand{\mV}{\mathcal{V}}
\newcommand{\menos}{\setminus}
\newcommand{\sone}{\mathsf{S}_1}
\newcommand{\gone}{\mathsf{G}_1}
\newcommand{\gfin}{\mathsf{G}_{fin}}
\newcommand{\sfin}{\mathsf{S}_\mathrm{fin}}
\newcommand{\0}{\mathrm{o}}
\newcommand{\mK}{\mathcal{K}}
\newtheorem*{claim}{Claim}
%%%%%
\newtheorem*{Prov}{Provision}
\newtheorem*{idfn}{Definition}
\newenvironment{definition}{\begin{idfn}
\rm}{\end{idfn}}

\section{Introduction}

Let $X$ be a topological space and let $\mathcal{O}$ be the set of all open covers for a space $X$. Given $\mC\in\mathcal{O}$ define $St(x,\mC)=\bigcup\{C\in\mC:x\in C\}$.

%Given a topological game $G(X)$, a strategy for player I is a function $\sigma$ definied over every finite sequence of moves of player II. 
%We denote by I $\uparrow G(X)$ if the player I has a winning strategy.

The \textbf{diagonal} of the $X\times X$ is the subset $\Delta=\{(x,x):x\in X\}$. We say that $X$ has the \textbf{$G_{\delta}$-diagonal property} if $\Delta$ is a $G_{\delta}$ subset of $X\times X$.

We say that $\mA\in\mathcal{O}$ is a \textbf{point-finite cover} if, for every $x\in X$, the set $\{A\in\mA:x\in A\}$ is finite. We say that a space $(X,\tau)$ is a \textbf{metacompact space} if every open cover has an open refinement that is point-finite. 

Let $X$ be a topological space and consider $L(X)=min\{\kappa\in\omega:$ given $\mC\in\mathcal{O}$ there is a $\mC'\subset\mC$ such that $\bigcup\mC'=X$ and $|\mC'|\leq\kappa\}$. We call this ordinal $L(X)$ the \textbf{Lindel\"of degree} of the space $X$.

Along this work we will use the standard topological definitions, following \cite{Engelking}. 

%In the second section we will use the Pressing Down Lemma to show that $\omega_1$ is not metacompact. And for this we use the notation of \cite{Kunen}.

%Let $\mC\subset\omega_1$.The set $\mC$ is a \textbf{cub} if it is closed and unbounded in $\omega_1$. We call $\mathcal{S}\subset\omega_1$ a \textbf{stationary set} if $\mathcal{S}\cap\mathcal{C}\neq\emptyset$ for every $\mathcal{C}$ cub's in $\omega_1$. Let $\kappa$ a regular cardinal and $f:S\rightarrow\kappa$. The function $f$ is a \textbf{regressive-function} if for every $\gamma\in S$ implies that $f(\gamma)<\gamma$.

%\begin{lem}[Pressing Down Lemma]\label{pressingdownlemma}\index{Lema do ``Pressing Down''}
%Let $\kappa>\omega$ regular cardinal, $S$ a stationary subset of $\kappa$ and $f:S\rightarrow\kappa$ a regressive function. Then there is $\alpha<\kappa$ such that $f^{-1}\{\alpha\}$ is a stationary set.
%\end{lem}

\section{Relations between $S_1(\mathcal{O},\mathcal{H})$ and $G_1(\mathcal{O},\mathcal{H})$}

Recall the following characterization for the $G_{\delta}$-diagonal property.

\

\begin{thm}[Ceder,\cite{Ceder}] Let $(X,\tau)$ be a topological space. Then $X$ has the $G_{\delta}$-diagonal property if, and only if, there is a countable sequence of open covers  $(\mC_n)_{n\in\omega}\subset\mO$ such that, for every $x,y\in X$, with $x\neq y$, there is a $k\in\omega$ such that $y\notin St(x,\mC_k)$. In other words, for each $x\in X$, $\bigcap_{n\in\omega} St(x,\mC_n)=\{x\}$.
%\begin{proof}

%You can check the proof in \cite{Ceder}.

%Suppose $\Delta$ is $G_{\delta}$, then for all $n\in\omega$ there are open sets $V_n\subset X^2$ such that $\Delta=\bigcap_{n\in\omega}V_n$. For each $n\in\omega$, let $\mC_n=\{A\in\tau:A\times A\subset V_n\}$. Note that $\mC_n\in\mO$. Let $x\in X$, then $(x,x)\in\Delta\subset V_n$, then there are open sets $U, V\in\tau$ with $(x,y)\in U\times V\subset V_n$ and this implies that $(x,y)\in (U\cap V)\times(U\cap V)\subset V_n$, therefore $x\in U\cap V\in\mC_n$. Now, we will show that $\{x\}=\bigcap_{n\in\omega} St(x,\mC_n)$. We have $\{x\}\subset\bigcap_{n\in\omega} St(x,\mC_n)$. Let $y\in\bigcap_{n\in\omega} St(x,\mC_n)$, for every $n\in\omega$ there is $A_n\in\mC_n$ such that $x,y\in A_n$ with implies that $(x,y)\in A_n\times A_n\subset V_n$ and then $(x,y)\in\bigcap_{n\in\omega}V_n=\Delta$. Thus $x=y$.

%Now, let $(\mC_n)_{n\in\omega}\subset\mO$ such that for each $x\in X$ we have $\{x\}=\bigcap_{n\in\omega} St(x,\mC_n)$. Put $V_n=\bigcup\{A\times A: A\in\mC_n\}$, so $\Delta\subset\bigcap_{n\in\omega}V_n$. Pick $(x,y)\in\bigcap_{n\in\omega}V_n$, then for all $n\in\omega$ there is $A_n\in\mC_n$ such that $(x,y)\in A_n\times A_n$ then $x,y\in A_n\in\mC_n$. But $\{x\}=\bigcap_{n\in\omega} St(x,\mC_n)$.
%\end{proof}
\end{thm}

\

This $G_{\delta}$-diagonal characterization give us motivation for a selection principle.

\

Let $(X,\tau)$ be a topological space. Let $\mH=\{R\in(\tau\menos\{\emptyset\})^{\omega}:|\bigcap R|\geq2\}$. The notation $S_1(\mathcal{O},\mathcal{H})$ abbreviates the following statement:

\

\noindent Given $(\mC_n)_{n\in\omega}\subset\mathcal{O}$, for each $n\in\omega$ there is a $C_n\in\mC_n$ such that $|\bigcap_{n\in\omega}C_n|\geq2$. 

\

%And this selective principle yields the game $G_1(\mathcal{O},\mathcal{H})$:
Associated to this principle, there is a game $G_1(\mathcal{O},\mathcal{H})$ defined as follows.

\

\noindent In the $n$-th inning, Player I plays a $\mC_n\in\mathcal{O}$ and Player II chooses a $C_n\in\mC_n$. At the end, Player II is the winner if $|\bigcap_{n\in\omega}C_n|\geq2$.

\

%$X^2$ has no $G_{\delta}$-diagonal is the same to say that: 

%\

%\noindent Given $(\mC_n)_{n\in\omega}$ there are $x,y\in X$ with $x\neq y$ such that for each $n\in\omega$ there is a $C_n\in\mC_n$ such that $x,y\in C_n$. Then $|\bigcap_{n\in\omega}C_n|\geq2$. And this is equal to $S_1(\mathcal{O},\mathcal{H})$.

\

Note that:

\

\begin{center}
$\Delta$ is $G_{\delta}$ in $X^2$ $\Leftrightarrow$ $\neg S_1(\mathcal{O},\mathcal{H})$ $\Rightarrow$ I $\uparrow G_1(\mathcal{O},\mathcal{H})$.
\end{center}

\

Where I $\uparrow G_1(\mathcal{O},\mathcal{H})$ means that Player I has a winning strategy for $G_1(\mathcal{O},\mathcal{H})$.

\

%Let's discuss what conditions we must have in $(X,\tau)$ to get the other implication.
In the following, we will discuss when the second implication can be reversed.

\

\begin{prop}\label{Lindelof}
Let $X$ be a Lindel\"of space. If Player I has a winning strategy in $G_1(\mathcal{O},\mathcal{H})$ then $S_1(\mathcal{O},\mathcal{H})$ does not hold.
\begin{proof}
Without loss of generality, we can assume that at each inning Player I plays a countable open covering. Then a winning strategy for Player I can be identified with a family $\{\mC_{\rho}:\rho\in\omega^{<\omega}\}$ such that $\mC_{\rho}=\{A_{\rho\smallfrown n}:n\in\omega\}$, where $\bigcup\{A_{\rho\smallfrown n}:n\in\omega\}=X$ and $|\bigcap\{A_{f\upharpoonright n}:n\in\omega\}|\leq1$ for every $f\in\omega^{\omega}$.

We will show that, for each $x\in X$, $\{x\}=\bigcap\{St(x,\mC_{\rho}):\rho\in\omega^{<\omega}\}$. Suppose that it does not happen. Then there are $x,y\in X$ with $x\neq y$ such that for every $\rho\in\omega^{<\omega}$, $y\in St(x,\mC_{\rho})$. Note that, the first move of Player I is $\{A_n:n\in\omega\}$. Then player II can choose $A_{n_{0}}$ such that $x,y\in A_{n_{0}}$. The next move of Player I is $\{A_{n_{0}n}:n\in\omega\}$ and Player II can choose $A_{n_{0}n_{1}}$ such that $x,y\in A_{n_{0}n_{1}}$ and so on. Then there is a $g\in\omega^{\omega}$ such that $x,y\in A_{g\upharpoonright n}$ for each $n\in\omega$, so $|\bigcap_{n\in\omega}A_{g\upharpoonright n}|\geq2$ which is a contradiction. %Therefore, there is a $\rho\in\omega^{<\omega}$ such that $y\notin St(x,\mC_{\rho})$.
\end{proof}
\end{prop}

%\begin{cor}
%Let $X$ be a Lindel\"of space then $|X|\leq2^{\omega}$. 
%\begin{proof}
%Let $x\in X$. In each inning player II chooses open set which contains the point $x$. By the Theorem \ref{Lindelof} the intersection of these open sets is the point $x$. 
%\end{proof}
%\end{cor}

\

The Proposition \ref{Lindelof} is not true for non Lindel\"of spaces as we will see in the following. Consider $\omega_1$ with the usual order topology. Then  %We will see below that in $\omega_1$ we have I $\uparrtrow G_1(\mathcal{O},\mathcal{H})$, but $\Delta$ is not $G_{\delta}$ in $\omega_1\times\omega_1$.

\

\begin{example}\label{ex}
Player I has a winning strategy in $G_1(\mathcal{O},\mathcal{H})$ played on $\omega_1$ and $S_1(\mathcal{O},\mathcal{H})$ holds.
\begin{proof}
Let $S$ be the set of all successors ordinals less than $\omega_1$ and $L=\omega_1\menos S$. For every $\gamma\in L$ pick a sequence of ordinals $\{\alpha_n^{\gamma}:n\in\omega\}\subset\omega_1$ such that $\alpha_n^{\gamma}<\alpha_{n+1}^{\gamma}<\gamma$ and $sup\{\alpha_n^{\gamma}:n\in\omega\}=\gamma$. For each $n\in\omega$ let $V_n^{\gamma}=]\alpha_n^{\gamma},\gamma]$ and $A=\{\{\alpha\}:\alpha\in S\}$. 

At the first inning Player I chooses $\mC_0=A\cup\{V_0^{\gamma}:\gamma\in L\}$. Note that if Player II chooses $\{\alpha\}$ for some $\alpha\in S$, then Player II loses the game. Then, we can suppose that Player II chooses $V_0^{\gamma_0}$ for some $\gamma_0\in L$. At the second inning Player I plays $\mC_1=A\cup\{]\gamma_0,\omega_1[\}\cup\{V_1^{\gamma}:\gamma\in L$ and $\gamma\leq\gamma_0\}$. Suppose Player II chooses $V_1^{\gamma_1}$ for some $\gamma_1\in L$ such that $\gamma_1\leq\gamma_0$. Note that, if $V_0^{\gamma_0}\cap V_1^{\gamma_1}=\emptyset$ then Player II loses. So, Player I chooses $\mC_2=A\cup\{]\gamma_1,\omega_1[\}\cup\{V_2^{\gamma}:\gamma\in L$ and $\gamma\leq\gamma_1\}$ and so on.

Since $\{\gamma_n:n\in\omega\}$ is a decreasing sequence of ordinals there is a $k\in\omega$ such that $\gamma_n=\gamma_k$ for every $n\geq k$. So, Player II chooses $V_n^{\gamma_k}\in\mC_n$ for every $n\geq k$. Therefore, $\bigcap_{n\in\omega}V_n^{\gamma_n}\subset\bigcap_{n\geq k}V_n^{\gamma_k}=\{\gamma_k\}$. % because $\{V_n^{\gamma_k}:n\in\omega\}$ is a local basis for $\gamma_k$.

\

Finally, note that $\omega_1$ is a countably compact non compact space, therefore, $\omega_1$ does not have the $G_{\delta}$-property, see e.g. \cite{Engelking}. Thus, $S_1(\mathcal{O},\mathcal{H})$ holds.

%and we have a result from \cite{Engelking} that says that if the space $X$ is countably compact and $\Delta$ is a $G_{\delta}$ in $X^2$ then $X$ is compact. 
\end{proof}
\end{example}

\

Note that the last proof works for every ordinal with uncountable cofinality. Therefore, the following is true.

\begin{prop}
Let $\alpha$ be an ordinal with uncountable cofinality. Then Player I has a winning strategy for $G_1(\mathcal{O},\mathcal{H})$ played on $A_{\alpha}=\{\beta<\alpha:cf(\beta)=\omega\}$. 
%\begin{proof}
%For each $\gamma\in A_{\alpha}$ pick a sequence of ordinals $\{\beta_n^{\gamma}:n\in\omega\}\subset A_{\alpha}$ such that $\beta_n^{\gamma}<\beta_{n+1}^{\gamma}<\gamma$ and $sup\{\beta_n^{\gamma}:n\in\omega\}=\gamma$. At the first inning Player I chooses $\mC_0=\{V_0^{\gamma}:\gamma\in A_{\alpha}\}$. Then, Player II chooses $V_0^{\gamma_0}$ for some $\gamma_0\in A_{\alpha}$. At the second inning Player I plays $\mC_1=\{]\gamma_0,\omega_1[\}\cup\{V_1^{\gamma}:\gamma\in A_{\alpha}$ and $\gamma\leq\gamma_0\}$. Suppose Player II choose $V_1^{\gamma_1}$ for some $\gamma_1\in A_{\alpha}$ such that $\gamma_1\leq\gamma_0$. Note that, if $V_0^{\gamma_0}\cap V_1^{\gamma_1}=\emptyset$ then II lose. So, Player I chooses $\mC_2=\{]\gamma_1,\omega_1[\}\cup\{V_2^{\gamma}:\gamma\in A_{\alpha}$ and $\gamma\leq\gamma_1\}$ and so on.

%We have a decreasing sequence of ordinals $\{\gamma_n:n\in\omega\}$, then there is a $k\in\omega$ such that $\gamma_n=\gamma_k$ for every $n\geq k$. So, Player II chooses $V_n^{\gamma_k}\in\mC_n$ for every $n\geq k$. Therefore, $\bigcap_{n\in\omega}V_n^{\gamma_n}\supset\bigcap_{n\geq k}V_n^{\gamma_k}=\{\gamma_k\}$, because $\{V_n^{\gamma_k}:n\in\omega\}$ is a local basis for $\gamma_k$.

%\end{proof}
\end{prop}

\

Now, we will see that the second implication can be reversed for hereditarily metacompact spaces. But, before that, we need some auxiliary results

\
%Besides, for each non countable ordinal $\alpha$ such that $cf(\alpha)>\omega$ we have that $\neg S_1(\mathcal{O},\mathcal{H})$ on $A_{\alpha}$.

%\section{Some results about Metacompact spaces}

Let $Y\subset X$ and let $\mathcal{O}(Y)$ be the set of all open covers for $Y$. Let $P(X)$ be the following game:
At the first inning Player I chooses $\mC_0\in\mathcal{O}$ and Player II answers by taking $A_0\in\mC_0$. At each inning $n\geq1$ Player I chooses $\mC_n\in\mathcal{O}(A_{n-1})$ and then Player II answers by taking $A_n\in\mC_n$. We say that Player II wins if $|\bigcap_{n\in\omega}A_n|\geq2$.

\

\begin{prop}
If Player I has a winning strategy in $G_1(\mathcal{O},\mathcal{H})$ then Player I has a winning strategy in $P(X)$ .
\begin{proof}
Let $\sigma$ be a winning strategy for Player I in $G_1(\mathcal{O},\mathcal{H})$. At the first inning of $P(X)$ Player I plays $\sigma(\emptyset)$ and Player II chooses $A_0\in\sigma(\emptyset)$. Then, at the second inning of $P(X)$, Player I plays $\mC_1=\{A_0\cap A:A\in\sigma(A_0)\}$ and Player II chooses $C_1\in\mC_1$.Note that $C_1=A_0\cap A_1$ for some $A_1\in\sigma(A_0)$. So, at the $n$th inning of the game $P(X)$, Player I plays $\mC_n=\{A_0\cap A_1\cap\cdots\cap A_{n-1}\cap C:C\in\sigma(A_0,...,A_{n-1})\}$ and Player II chooses $C_n\in\mC_n$, with $C_n=A_0\cap A_1\cap\cdots\cap A_n$ for some $A_n\in\sigma(A_0\cap A_1\cap\cdots\cap A_{n-1})$. Therefore $|\bigcap_{n\in\omega}C_n|=|\bigcap_{n\in\omega}A_n|\leq1$, where $A_0=C_0$.
\end{proof}
\end{prop}

\

\begin{lem}
Let $X$ be a hereditarily metacompact space. If Player I has a winning strategy in $P(X)$ then there is a winning strategy for Player I in $P(X)$ such that Player I only plays point-finite open covers.
\begin{proof}
Let $\sigma$ be a winning strategy for Player I in $P(X)$. Let $\sigma(\emptyset)$ be the first move of Player I and let $\mC_0=\sigma^*(\emptyset)$ be a point-finite refinement of $\sigma(\emptyset)$. If Player II chooses $A_0^*\in\mC_0$, then there is an $A_0\in\sigma(\emptyset)$ such that $A_0^*\subset A_0$. Let $\sigma^*(A_0)$ be a point-finite refinement of $\sigma(A_0)$. Let $\mC_1=\{B^*\cap A_0^*:B^*\in\sigma^*(A_0)\}$ be the play for Player I. Note that $\mC_1$ is a point-finite cover for $A_0^*$. Then Player II chooses $A_1^*\in\mC_1$ such that $A_1^*=B^*_1\cap A_0^*$ with $B^*_1\in\sigma^*(A_0)$, then there is an $A_1\in\sigma(A_0)$ such that $B^*_1\subset A_1$. Let $\sigma^*(A_0,A_1)$ be a point-finite refinement of $\sigma(A_0,A_1)$ and Player I plays $\mC_2=\{B^*\cap A_1^*:B^*\in\sigma^*(A_0,A_1)\}$. Again, note that $\mC_2$ is a point-finite cover for $A_1^*$ . Then Player II chooses $A^*_2\in\mC_2$, $A^*_2=B^*_2\cap A_1^*$ with $B_2^*\in\sigma^*(A_0,A_1)$, then there is an $A_2\in\sigma(A_0,A_1)$ such that $B^*\subset A_2$.

Proceeding this way, in the $n$-th inning Player I plays $\mC_n=\{B^*\cap A_{n-1}^*:B^*\in\sigma^*(A_0,...,A_{n-1})\}$ and Player II chooses $A_n^*\in\mC_n$. Note that for each $n\in\omega$ $A_n^*\subset A_n$. Therefore,
\[|\bigcap_{n\in\omega}A_n^*|\leq|\bigcap_{n\in\omega}A_n|\leq1\]
\end{proof}
\end{lem}

\begin{prop}
If $(X,\tau)$ is a hereditarily metacompact space and Player I has a winning strategy in $P(X)$ then $S_1(\mathcal{O},\mathcal{H})$ does not hold.
\begin{proof}
As we saw above we can suppose that all covers played by Player I are point-finite. Let $\sigma$ be a winning strategy for Player I in $P(X)$. For each $x\in X$ let $S(x,\emptyset)=\{C\in\sigma(\emptyset):x\in C\}$ and $S(x,C_0,...,C_n)=\{C\in\sigma(C_0,...,C_n):C_n\in S(x,C_0,...,C_{n-1}),...,C_0\in S(x,\emptyset)$ and $x\in C\}$. Note that $S(x,\emptyset)$ and $S(x,C_0,...,C_n)$ are finite sets. For every $x\in X$, let $V_0^x=\bigcap S(x,\emptyset)$ and $V_n^x=\bigcap\{\bigcap S(x,C_0,...,C_n):C_n\in S(x,C_0,...,C_{n-1}),...,C_0\in S(x,\emptyset)\}$.

For each $n\in\omega$ we define $\mC_n=\{V_n^z:z\in X\}$. We will show that $\bigcap_{n\in\omega}St(x,\mC_n)=\{x\}$. Suppose that it does not happen, then there are $x,y\in X$ with $x\neq y$ such that $y\in\bigcap_{n\in\omega}St(x,\mC_n)$. Note that if $St(x,\mC_n)=\bigcup\{V_n^z:z\in X $ and $x\in V_n^z\}$ then for every $n\in\omega$ there is $z_n\in X$ such that $x,y\in V_n^{z_n}$. Let $Lev(n)=\bigcup\{S(x,C_0,...,C_n):C_n\in S(x,C_0,...,C_{n-1}),...,C_0\in S(x,\emptyset)\}$ and let $T=\bigcup_{n\in\omega}Lev(n)$. Note that $(T,\leq)$, is a tree ordered by ``$\supseteq$''. Note that $V_n^{z_n}\in Lev(n)$ for each $n\in\omega$.

\begin{claim}
There is a branch $R$ of $(T,\leq)$ such that $x,y\in\bigcap R$.
\begin{proof}
Every level of the tree $(T,\leq)$ has finitely many elements and each element of a level forks in another finitely many elements of the next level. Then there are $C_0\in Lev(0)$ and $\mA_0\subset T$ with $|\mA_0|=\omega$ such that, for every $A\in\mA_0$, $C_0\leq A$. There are $C_1\in Lev(1)$ and $\mA_1\subset\mA_0$ with $|\mA_1|=\omega$ such that, for every $A\in\mA_1$, $C_1\leq A$ and $C_0\leq C_1$. Proceeding this way, we can find for every $n$ a $C_n\in Lev(n)$ such that $C_0\leq C_1\leq\cdots\leq C_n$ and an $\mathcal{A}_n\subset\mathcal{A}_{n-1}$ such that $|\mathcal{A}_n|=\omega$ and, for each $A\in\mathcal{A}_n$, $C_n\leq A$. So, $R=(C_n)_{n\in\omega}$ is a branch and $x,y\in C_n$ for every $n\in\omega$.
\end{proof}
\end{claim}

Therefore, Player II wins, which is a contradiction. Then, $\bigcap_{n\in\omega}St(x,\mC_n)=\{x\}$ for every $x\in X$.

\end{proof}
\end{prop}

\begin{cor}
If $(X,\tau)$ is a hereditarily metacompact space and Player I has a winning strategy in $G_1(\mathcal{O},\mathcal{H})$ then $S_1(\mathcal{O},\mathcal{H})$ does not hold.
\end{cor}

Note that the Pixley-Roy hyperspaces are always hereditarily metacompact, so we have the following corollary:

\begin{cor}
In any Pixley-Roy hyperspace, $S_1(\mathcal{O},\mathcal{H})$ holds if, and only if, Player I does not have a winning strategy in $G_1(\mathcal{O},\mathcal{H})$.
\end{cor}

\section{Applications}

%In this section we will see some applications of the game $G_1(\mathcal{O},\mathcal{H})$.

In \cite{Buzyakova} it is shown that every space with the countable chain condition and the regular $G_{\delta}$-diagonal property has size at most $\frak{c}$. We will show a similar result involving the game $G_1(\mathcal{O},\mathcal{H})$. After that, we will see others applications.

\

Let $(X,\tau)$ be a topological space. Let $\mathcal{S}=\{(S_n)_{n\in\omega}:S_n\in\tau\menos\{\emptyset\}\}$. In the following we will denote each sequence $(S_n)_{n\in\omega}$ only by $S$. For each $A,B\in\mathcal{S}$ with $A\neq B$ let $d:\mathcal{S}\times\mathcal{S}\rightarrow[0,1]$ be a function such that $d(A,B)=1/(n+1)$ where $n=min\{k\in\omega:A_k\neq B_k\}$ and $d(A,A)=0$. Note that $(\mathcal{S},d)$ is a metric space.

\

Consider $G_1^*(\mathcal{O},\mathcal{H})$ the following game: In the $n$-th inning, Player I plays a $\mC_n\in\mathcal{O}$ and Player II chooses a $C_n\in\mC_n$. At the end, Player II is the winner if there is a $k\in\omega$ such that $|\bigcap_{k\leq n}C_n|\geq2$.

Observe that if Player I has a winning strategy in $G_1(\mathcal{O},\mathcal{H})$ then there is a winning strategy for Player I in $G_1^*(\mathcal{O},\mathcal{H})$. Indeed, let $\sigma$ be a winning strategy for Player I in $G_1(\mathcal{O},\mathcal{H})$. Let us define a strategy for Player I in the $G_1^*(\mathcal{O},\mathcal{H})$ in the following way. At the $n$-th inning, Player I chooses:
\begin{itemize}
\item $\sigma(\emptyset)$, if $n=p$ for $p$ a prime number.
\item $\sigma(C_{p},...,C_{p^{k-1}})$, if $n=p^k$ for $p$ a prime number and for some $k\in\omega$, $k>1$.
\item $\sigma(\emptyset)$, if $n$ is not a power of a prime number.
\end{itemize}
%Note that for each prime number $p\geq2$ we have a different game $G_1(\mathcal{O},\mathcal{H})$ where for each inning is a power of $p$ and Player I has a winning strategy.  

 Note that, at the $n$-th inning, Player II chooses $C_n$. So, for any $k\in\omega$ there is a prime number such that $p>k$, then
\[|\bigcap_{n\geq k}C_n|\leq|\bigcap_{n\geq p}C_n|\leq|\bigcap_{n\geq1}C_{p^n}|\leq1.\]

\

\begin{prop}
Let $(X,\tau)$ be a metacompact and separable space. If Player I has a winning strategy in $G_1(\mathcal{O},\mathcal{H})$ then $|X|\leq\frak{c}$.
\begin{proof}
Let $\sigma$ be the winning strategy for Player I. For each $x\in X$ let $S(x,\emptyset)=\{C\in\sigma(\emptyset):x\in C\}$ and $S(x,C_0,...,C_{n-1})=\{C_n\in\sigma(C_0,...,C_{n-1}):C_{n-1}\in S(x,C_0,...,C_{n-2}),...,C_0\in S(x,\emptyset)$, and $x\in C_n\}$. For each $x\in X$ let $\mathcal{S}_x=\{R\in\mathcal{S}:R_n\in\sigma(R_0,...,R_{n-1}),...,R_0\in\sigma(\emptyset)$ and $\bigcap R=\{x\}\}$. For each $x\in X$ and for each $k\in\omega$ let $\mathcal{U}_x^k=\{R\in\mathcal{S}:R_n\in S(x,C_0,...,C_{n-1}),...,C_0\in S(x,\emptyset)$ and $R_n=A_n$ for every $n\geq k$ and for some $A\in\mathcal{S}_x\}$. Let $\mathcal{U}_x=\bigcup_{n\in\omega}\mathcal{U}_x^k$ and let $\mathcal{U}=\bigcup_{x\in X}\mathcal{U}_x$. Note that $(\mathcal{U},d\upharpoonright\mathcal{U})$ is a metric space and there is a sobrejective function $g:\mathcal{U}\rightarrow X$ such that $g(U)=\bigcap U$.

Let $D$ be a countable and dense subset of $X$. For each $d\in D$ fix $A^d\in\mathcal{S}_d$ and for each $n\in\omega$ let $\mathcal{D}_d^n=\{H\in\mathcal{U}_d:H_k=A^d_k$ for each $k\geq n\}$ and let $\mathcal{D}_d=\bigcup_{n\in\omega}\mathcal{D}_d^n$. So, define $\mathcal{D}=\bigcup_{d\in D}\mathcal{D}_d$. Observe that $\mathcal{D}$ is a countable subset of $\mathcal{U}$. We will show that $\overline{\mathcal{D}}=\mathcal{U}$. Let $Y\in\mathcal{U}_y$ we will show that $B(Y,1/(n+1))\cap\mathcal{D}\neq\emptyset$, where $B(Y,1/(n+1))$ is the open ball of center $Y$ and radius $1/(n+1)$. For each $x\in X$ and $k\in\omega$ let $V_k^x=\bigcap\{\bigcap S(x,C_0,...,C_{k-1}):C_{k-1}\in S(x,C_0,...,C_{k-1}),...,C_0\in S(x,\emptyset)\}$. Since $X$ is a metacompact space, $\bigcap_{k\leq n}V_k^y$ is a non empty open set then let $d\in (\bigcap_{k\leq n}V_k^y)\cap D$. So $S(d,C_0,...,C_k)\supset S(y,C_0,...,C_k)$ for each $k\leq n$. Let $H_k=Y_k$ for $k\leq n$ and let $H_k=A^d_k$ for $k>n$. Therefore $H\in B(Y,1/(n+1))\cap\mathcal{D}$. Since $\mathcal{U}$ is a metric space, and a separable then $\mathcal{U}$ has a countable base. So, $|\mathcal{U}|\leq\frak{c}$ and since $g:\mathcal{U}\rightarrow X$ is a sobrejective function the result follows.
\end{proof}
\end{prop}

Note that the metacompactness is important. Consider the Katetov's extension $K(\omega)$ of $\omega$. This is a separable space with $G_{\delta}$-diagonal and have cardinality bigger than $\frak{c}$. %So we conclude that $K(\omega)$ is not a metacompact space.

\

\begin{prop}
Let $\kappa$ be the Lindel\"of degree of a space $X$. If Player I has a winning strategy for $G_1(\mathcal{O},\mathcal{H})$  then $|X|\leq\kappa^{\omega}$.
\begin{proof}
Let $x\in X$. At each inning, we can suppose that Player I plays a cover of size $\kappa$ and Player II chooses an open set which contains the point $x$. Then, since Player I has a winning strategy in $G_1(\mathcal{O},\mathcal{H})$, the intersection of all these open sets is $\{x\}$. The size of the set of all branches provided by the winning strategy of Player I is at most $\kappa^{\omega}$.
\end{proof}
\end{prop}

\

%As we saw at the Proposition \ref{Lindelof} the lenght of the game and the Lindelöf degree are essential to space $X$ has a $G_{\delta}$-diagonal property. 

Let $\kappa$ be an infinite cardinal. Consider the selection principle $S_1^{\kappa}(\mathcal{O},\mathcal{H})$ given by the following statement:

\

Given $(\mC_{\xi})_{\xi<\kappa}\subset\mathcal{O}$, for each $\xi<\kappa$ there is a $C_{\xi}\in\mC_{\xi}$ such that $|\bigcap_{\xi<\kappa}C_{\xi}|\geq2$. 

\

Again, we have associated to this principle a game $G_1^{\kappa}(\mathcal{O},\mathcal{H})$ defined as follows.

\

\noindent In the $\xi$-th inning, Player I plays a $\mC_{\xi}\in\mathcal{O}$ and Player II chooses a $C_{\xi}\in\mC_{\xi}$. At the end, Player II is the winner if $|\bigcap_{\xi<\kappa}C_{\xi}|\geq2$.

\

Note that:

\

\begin{center}
$\Delta$ is $G_{\kappa}$ in $X^2$ $\Leftrightarrow$ $\neg S_1^{\kappa}(\mathcal{O},\mathcal{H})$ $\Rightarrow$ I $\uparrow G_1^{\kappa}(\mathcal{O},\mathcal{H})$.
\end{center}

\

In the following we show that a winning strategy for Player I in $G_1(\mathcal{O},\mathcal{H})$ gives a bound for which selections of the form $S_1^{\kappa}(\mathcal{O},\mathcal{H})$ can hold.

\

\begin{prop}
Let $\kappa$ be the Lindel\"of degree of a space $X$. If Player I has a winning strategy in $G_1(\mathcal{O},\mathcal{H})$ then $S_1^{\kappa}(\mathcal{O},\mathcal{H})$ does not hold.
\begin{proof}
Without loss of generality, we can assume that at each inning Player I plays an open covering of size $\kappa$. Consider a winning strategy for Player I given by the $\{\mC_{\rho}:\rho\in\kappa^{<\omega}\}$ such that $\mC_{\rho}=\{A_{\rho^{\smallfrown}\xi}:\xi<\kappa\}$, $\bigcup\{A_{\rho^{\smallfrown}\xi}:\xi<\kappa\}=X$ and $|\bigcap\{A_{f\upharpoonright\xi}:\xi<\kappa\}|\leq1$ for every $f\in\kappa^{\omega}$.

Note that $|\kappa^{< \omega}|=\kappa$. Suppose that $S_1^{\kappa}(\mathcal{O},\mathcal{H})$ holds, i.e., there are $x,y\in X$ with $x\neq y$ such that for every $\rho\in\kappa^{<\omega}$ we have $y\in St(x,\mC_{\rho})$. Then there is a $g\in\kappa^{\omega}$ such that $x,y\in A_{g\upharpoonright\xi}$ for each $\xi<\kappa$, so $|\bigcap_{\xi<\kappa}A_{g\upharpoonright\xi}|\geq2$. Therefore, the Player I does not have winning strategy in $G_1(\mathcal{O},\mathcal{H})$.
\end{proof}
\end{prop}

%Note that the last result holds to every non countable ordinal with countable cofinality.
%We can extend this last result for sets which their cardinals numbers as big as we want. In other words, given a cardinal $\kappa$ we can get a set of size bigger than $\kappa$ such that Player I has a winning strategy in $G_1(\mathcal{O},\mathcal{H})$ and $S_1(\mathcal{O},\mathcal{H})$.

%\begin{remark}
%With the last result we can show that $\omega_1$ is not metacompact.
%\end{remark}

\section*{Acknowledgement}

We would like to thank Professor Angelo Bella who gave to us several useful suggestions. And we also thank the anonymous referee for the careful corrections and useful suggestions.


\begin{thebibliography}{9}

\bibitem{Buzyakova} R.~Buzyakova. \emph{Cardinalities of ccc-spaces with regular Gδ-diagonals}. Topology Appl. 153 (2006), no. 11, 1696–1698.

\bibitem{Ceder} J. G.~Ceder. \emph{Some generalizations of metric spaces}. Pacific J. Math., 11 1961, 105-125.

\bibitem{Engelking} R.~Engelking. \emph{General topology}. Heldermann Verlag, 1989.

\bibitem{Kunen} K.~Kunen. \emph{Set Theory: An Introduction to Independence Proofs}. North-Holland Publishing Company, 1980.

%\bibitem{Zenor} P.~Zenor. \emph{Spaces with regular $G_{\delta}$-diagonals. General Topology and its
%Relations to Modern Analysis and Algebra III}, (Proc. Third Prague Topological Sympos., 1971), pp. 471-473. Academia, Prague, 1972. 

%\bibitem{Arhan1} A.~Arhangel'skii. \emph{The frequency spectrum of a topological space and the product operation}. Trans. Moscow Math. Soc. 1981, 2: 163-200.

%\bibitem{Aurichi} L.~F.~Aurichi, A.~Bella, A. \emph{Productively countably tight spaces and topological games.} arXiv:1307.7928, 2013.


%\bibitem{Grif} D.~F.~Griffiths, D.~J.~Higham,
%\emph{learning \LaTeX, }Siam, 1997.
 %\bibitem{Lam}
%L.~Lamport, \LaTeX: \emph{ A Document Preparation
%System. User's Guide and Reference Manual, }2nd
%edition, Addison-Wesley, Reading, MA, 1994.

\end{thebibliography}
 \end{document}